\documentclass[a4]{article}
\usepackage{amssymb}
\usepackage{diagrams}
\begin{document}
\newcommand{\nc}[2]{\newcommand{#1}{#2}}
\newcommand{\rnc}[2]{\renewcommand{#1}{#2}}
\rnc{\theequation}{\thesection.\arabic{equation}}
\def\note#1{{}}
\def\Label#1{\label{#1}\ifmmode\llap{[#1] }\else
\marginpar{\smash{\hbox{[#1]}}}\fi}

\newtheorem{de}{Definition $\!\!$}[section]
\newtheorem{pr}[de]{Proposition $\!\!$}
\newtheorem{lem}[de]{Lemma $\!\!$}
\newtheorem{corollary}[de]{Corollary $\!\!$}
\newtheorem{them}[de]{Theorem $\!\!$}
\newtheorem{example}[de]{\it Example $\!\!$}
\newtheorem{remark}[de]{\it Remark $\!\!$}

\nc{\beq}{\begin{equation}}
\nc{\eeq}{\end{equation}}
\rnc{\[}{\beq}
\rnc{\]}{\eeq}
\rnc{\Label}{\label}
\nc{\ba}{\begin{array}}
\nc{\ea}{\end{array}}
\nc{\bea}{\begin{eqnarray}}
\nc{\beas}{\begin{eqnarray*}}
\nc{\eeas}{\end{eqnarray*}}
\nc{\eea}{\end{eqnarray}}
\nc{\be}{\begin{enumerate}}
\nc{\ee}{\end{enumerate}}
\nc{\bd}{\begin{diagram}}
\nc{\ed}{\end{diagram}}
\nc{\bi}{\begin{itemize}}
\nc{\ei}{\end{itemize}}
\nc{\bpr}{\begin{pr}}
\nc{\bth}{\begin{them}}
\nc{\ble}{\begin{lem}}
\nc{\bco}{\begin{corollary}}
\nc{\bre}{\begin{remark}}
\nc{\bex}{\begin{example}}
\nc{\bde}{\begin{de}}
\nc{\ede}{\end{de}}
\nc{\epr}{\end{pr}}
\nc{\ethe}{\end{them}}
\nc{\ele}{\end{lem}}
\nc{\eco}{\end{corollary}}
\nc{\ere}{\hfill\mbox{$\losenge$}\end{remark}}
\nc{\eex}{\hfill\mbox{$\losenge$}\end{example}}
\nc{\bpf}{{\it Proof.~~}}
\nc{\epf}{\hfill\mbox{$\square$}\vspace*{3mm}}
\nc{\hsp}{\hspace*}
\nc{\vsp}{\vspace*}

\nc{\ot}{\otimes}
\nc{\te}{\!\ot\!}
\nc{\bmlp}{\mbox{\boldmath$\left(\right.$}}
\nc{\bmrp}{\mbox{\boldmath$\left.\right)$}}
\nc{\LAblp}{\mbox{\LARGE\boldmath$($}}
\nc{\LAbrp}{\mbox{\LARGE\boldmath$)$}}
\nc{\Lblp}{\mbox{\Large\boldmath$($}}
\nc{\Lbrp}{\mbox{\Large\boldmath$)$}}
\nc{\lblp}{\mbox{\large\boldmath$($}}
\nc{\lbrp}{\mbox{\large\boldmath$)$}}
\nc{\blp}{\mbox{\boldmath$($}}
\nc{\brp}{\mbox{\boldmath$)$}}
\nc{\LAlp}{\mbox{\LARGE $($}}
\nc{\LArp}{\mbox{\LARGE $)$}}
\nc{\Llp}{\mbox{\Large $($}}
\nc{\Lrp}{\mbox{\Large $)$}}
\nc{\llp}{\mbox{\large $($}}
\nc{\lrp}{\mbox{\large $)$}}
\nc{\lbc}{\mbox{\Large\boldmath$,$}}
\nc{\lc}{\mbox{\Large$,$}}
\nc{\Lall}{\mbox{\Large$\forall\;$}}
\nc{\bc}{\mbox{\boldmath$,$}}
\nc{\ra}{\rightarrow}
\nc{\ci}{\circ}
\nc{\cc}{\!\ci\!}
\nc{\lra}{\longrightarrow}
\rnc{\to}{\mapsto}
\nc{\imp}{\Rightarrow}
\rnc{\iff}{\Leftrightarrow}
\nc{\inc}{\mbox{$\,\subseteq\;$}}
\rnc{\subset}{\inc}
\nc{\0}{\sp{(0)}}
\nc{\1}{\sp{(1)}}
\nc{\2}{\sp{(2)}}
\nc{\3}{\sp{(3)}}
\nc{\4}{\sp{(4)}}
\nc{\5}{\sp{(5)}}
\nc{\6}{\sp{(6)}}
\nc{\7}{\sp{(7)}}
\newcommand{\Boxneu}{\square}
\def\tr{{\rm tr}}
\def\Tr{{\rm Tr}}
\def\st{\stackrel}
\def\<{\langle}
\def\>{\rangle}
\def\d{\mbox{$\mathop{\mbox{\rm d}}$}}
\def\id{\mbox{$\mathop{\mbox{\rm id}}$}}
\def\ker{\mbox{$\mathop{\mbox{\rm Ker$\,$}}$}}
\def\hom{\mbox{$\mathop{\mbox{\rm Hom}}$}}
\def\corep{\mbox{$\mathop{\mbox{\rm Corep}}$}}
\def\im{\mbox{$\mathop{\mbox{\rm Im}}$}}
\def\map{\mbox{$\mathop{\mbox{\rm Map}}$}}
\def\o{\sp{[1]}}
\def\t{\sp{[2]}}
\def\mo{\sp{[-1]}}
\def\z{\sp{[0]}}

\nc{\spp}{\mbox{${\cal S}{\cal P}(P)$}}
\nc{\ob}{\mbox{$\Omega\sp{1}\! (\! B)$}}
\nc{\op}{\mbox{$\Omega\sp{1}\! (\! P)$}}
\nc{\oa}{\mbox{$\Omega\sp{1}\! (\! A)$}}
\nc{\dr}{\mbox{$\Delta_{R}$}}
\nc{\dsr}{\mbox{$\Delta_{\Omega^1P}$}}
\nc{\ad}{\mbox{$\mathop{\mbox{\rm Ad}}_R$}}
\nc{\as}{\mbox{$A(S^3\sb s)$}}
\nc{\bs}{\mbox{$A(S^2\sb s)$}}
\nc{\slc}{\mbox{$A(SL(2,\C))$}}
\nc{\suq}{\mbox{$\cO(SU_q(2))$}}
\nc{\tc}{\widetilde{can}}
\def\slq{\mbox{$\cO(SL_q(2))$}}
\def\asq{\mbox{$\cO(S_{q,s}^2)$}}
\def\esl{{\mbox{$E\sb{\frak s\frak l (2,{\Bbb C})}$}}}
\def\esu{{\mbox{$E\sb{\frak s\frak u(2)}$}}}
\def\ox{{\mbox{$\Omega\sp 1\sb{\frak M}X$}}}
\def\oxh{{\mbox{$\Omega\sp 1\sb{\frak M-hor}X$}}}
\def\oxs{{\mbox{$\Omega\sp 1\sb{\frak M-shor}X$}}}
\def\Fr{\mbox{Fr}}

\rnc{\epsilon}{\varepsilon}
\rnc{\phi}{\varphi}
\nc{\ha}{\mbox{$\alpha$}}
\nc{\hb}{\mbox{$\beta$}}
\nc{\hg}{\mbox{$\gamma$}}
\nc{\hd}{\mbox{$\delta$}}
\nc{\he}{\mbox{$\varepsilon$}}
\nc{\hz}{\mbox{$\zeta$}}
\nc{\hs}{\mbox{$\sigma$}}
\nc{\hk}{\mbox{$\kappa$}}
\nc{\hm}{\mbox{$\mu$}}
\nc{\hn}{\mbox{$\nu$}}
\nc{\hl}{\mbox{$\lambda$}}
\nc{\hG}{\mbox{$\Gamma$}}
\nc{\hD}{\mbox{$\Delta$}}
\nc{\hT}{\mbox{$\Theta$}}
\nc{\ho}{\mbox{$\omega$}}
\nc{\hO}{\mbox{$\Omega$}}
\nc{\hp}{\mbox{$\pi$}}
\nc{\hP}{\mbox{$\Pi$}}

\nc{\qpb}{quantum principal bundle}
\def\gal{-Galois extension}
\def\hge{Hopf-Galois extension}
\def\ses{short exact sequence}
\def\csa{$C^*$-algebra}
\def\ncg{noncommutative geometry}
\def\wrt{with respect to}
\def\Ha{Hopf algebra}

\def\C{{\Bbb C}}
\def\N{{\Bbb N}}
\def\R{{\Bbb R}}
\def\Z{{\Bbb Z}}
\def\T{{\Bbb T}}
\def\cO{{\mathcal O}}
\def\O{\cO}
\def\cT{{\cal T}}
\def\cK{{\cal K}}
\def\cH{{\cal H}}
\def\ch{{\cal H}}
\def\H{{\cal H}}
\def\ta{\tilde a}
\def\tb{\tilde b}
\def\td{\tilde d}
\newcommand{\weg}[1]{}
{\noindent\LARGE\bf
Chern numbers for two families of noncommutative \\[.2cm]Hopf fibrations
\footnote{See the home page
of P.M.H.\ for a preliminary detailed version of this work.}}\\ \medskip

\noindent{\large\bf
Piotr M.~Hajac$^a$, Rainer Matthes$^b$, Wojciech Szyma\'nski$^c$\\
\medskip\medskip

\footnotesize\bf\begin{itemize}
\vspace*{-4mm}\item[$^a$]
Mathematisches Institut, Universit\"at M\"unchen,
Theresienstr.\ 39, M\"unchen, 80333, Germany;
and
Instytut Matematyczny, Polska Akademia Nauk,
ul.\ \'Sniadeckich 8, Warszawa, 00-950 Poland;
and
Katedra Metod Matematycznych Fizyki, Uniwersytet Warszawski,
ul.\ Ho\.za 74, Warszawa, 00-682 Poland;
URL: http://www.fuw.edu.pl/$\!\widetilde{\phantom{m}}\!$pmh\\

\vspace*{-4mm}\item[$^b$]
Fachbereich 2, TU Clausthal,
Leibnizstr.\ 4, D-38678 Clausthal-Zellerfeld, Germany;
E-mail: ptrm@pt.tu-clausthal.de\\

\vspace*{-4mm}\item[$^c$]
School of Mathematical and Physical Sciences, The University of Newcastle,
Callaghan, NSW 2308, Australia;
E-mail: wojciech@frey.newcastle.edu.au
\end{itemize}
}
\bigskip

\ $\!\!$\hrulefill\
\medskip

\noindent{\bf Abstract.}

{\small
We consider noncommutative line bundles associated with the Hopf fibrations of
$SU_q(2)$ over all Podle\'s spheres and with a locally trivial Hopf fibration
of $S^3_{pq}$. These bundles are given as finitely generated projective
modules associated via 1-dimensional representations of $U(1)$
with Galois-type extensions encoding the  principal fibrations of $SU_q(2)$
and $S^3_{pq}$.
We show that the Chern numbers of these modules coincide with the winding
numbers of representations defining them.
}\medskip\ ~\\

\weg{
{\noindent\large\bf
Chern numbers for two families of noncommutative Hopf fibrations}\\

\noindent{\bf R\'esum\'e.}

{\small
We consider noncommutative line bundles associated with the Hopf fibrations of
$SU_q(2)$ over all Podle\'s spheres and with a locally trivial Hopf fibration
of $S^3_{pq}$. These bundles are given as finitely generated projective
modules associated via 1-dimensional representations of $U(1)$
with Galois-type extensions encoding the  principal fibrations of $SU_q(2)$
and $S^3_{pq}$.
We show that the Chern numbers of these modules coincide with the winding
numbers of representations defining them.
}

\ $\!\!$\hrulefill\
\bigskip

\noindent{\large\bf Version fran\c{c}aise abr\'eg\'ee}\\
}
\noindent{\large\bf  Abridged version    }\\

In this paper, we combine the algebraic tool of Galois-type extensions
with the analytic tool of the noncommutative index formula to
study two kinds of examples of quantum fibrations.
Our main result is  that the line bundles associated to these principal
fibrations are mutually non-isomorphic.  This gives an estimate of
the positive cones of the algebraic $K_0$-groups of the base-space
quantum spheres.

Let $B\inc P$ be an inclusion of algebras such that $B$ is the coinvariant
subalgebra for some coalgebra $C$ coaction $\dr:P\ra P\ot C$.  Using the
framework of Galois-type extensions, one can say when such an extension of
algebras is principal. (The definition is tuned in such a way that for
commutative algebras it coincides with the concept of affine group scheme
torsors --- the principal bundles of algebraic geometry.) Every
principal $C$-extension $B\inc P$
allows one to assign to any finite-dimensional
corepresentation  of $C$ a finitely generated projective left $B$-module
of colinear homomorphisms    $\hom^C(V_\phi,P)$ \cite{bh}.
Taking its class in
$K_0(B)$ and composing it with the Chern character defines the
Chern-Galois character  from
the space of all finite-dimensional corepresentations of
$C$ to the even cyclic homology of $B$. On the other hand,
the $K$-homology Chern character assigns to finitely summable Fredholm modules
cyclic cocycles. In the 1-summable case it takes a particularly simple form,
notably it turns a pair of bounded $*$-representations
$(\rho_1,\rho_2)$ into a trace (cyclic 0-cocycle) on $B$ via
the formula $\tr_\rho=\Tr\ci(\rho_1-\rho_2)$.
The evaluation of this trace on  the Chern-Galois character
applied to a corepresentation gives a numerical invariant of the
$K_0$-class of the module  defined by this corepresentation.
Moreover, for our examples, the integrality of these invariants
(guaranteed by the noncommutative index formula) makes them computable.

Our first example of a noncommutative Hopf fibration has its source
in Klimek-Lesniewski quantum discs and the idea of local triviality.
Here we have a principal $\O(U(1))$-extension
$\O(S^2_{pq})\inc \O(S^3_{pq}))$ \cite{hms2},
  and every 1-dimensional corepresentation
$\phi_\mu(1)=1\ot u^{-\mu}$ yields a
finitely generated projective $\O(S^2_{pq})$-module (quantum line bundle)
$\O(S^3_{pq})_\mu:=\hom^{\O(U(1))}(\C_{\phi_\mu},\O(S^3_{pq}))$.
On the other hand, we have the following irreducible bounded
$*$-representations of $\cO(S^2_{pq})$
\cite{cm00}:
\bea
 &&\label{ro1}
\rho_1(f_0)e_k=(1-p^k)e_k,\;\;\;
\rho_1(f_1)e_k=\sqrt{1-p^{k+1}}\:e_{k+1},\;\;\;k\geq 0;
\\ &&\label{ro2}
\rho_2(f_0)e_k=e_k,\;\;\;
\rho_2(f_1)e_k=\sqrt{1-q^{k+1}}\:e_{k+1},\;\;\;k\geq 0.
\eea
Here $\{e_k\}_{k\geq 0}$ is an orthonormal basis of a separable Hilbert space
and $f_0,f_1$ are generators of the $*$-algebra $\cO(S^2_{pq})$.  These
representations give a desired trace, and
bring us to our first main result:
\bth
For all $\mu\in\Z$, the pairing between the cyclic 0-cocycle
${\rm tr}_\rho$ and the $K_0$-class of $\O(S^3_{pq})_\mu$ (Chern number)
coincides with the
winding number $\mu$, i.e.,
$
\langle \tr_\rho,[\O(S^3_{pq})_{\mu}]\rangle = \mu.
$
\ethe

Our second example is a family of noncommutative Hopf fibrations of $SU_q(2)$
over all Podle\'s quantum spheres $S^2_{q,s}$, $s\in[0,1]$.
As handling the generic Podle\'s spheres
requires going beyond the Hopf-Galois theory, they were among main motivating
examples driving the development of the theory of
principal coalgebra extensions of
noncommutative rings.
Now, much as before, we associate to every corepresentation $\phi_\mu$
a finitely generated projective left $\O(S^2_{q,s})$-module
$\O(SU_q(2))_{\mu,s}:=\hom^{\cO(SU_q(2))/J_s}(\C_{\phi_\mu},\cO(SU_q(2)))$.
Here    $J_s$ is
the coideal right ideal generated by $K$, $L-s$, $L^*-s$, and $K,L$
generate the polynomial algebras  $\O(S^2_{q,s})$.
One can show
that $\cO(SU_q(2))/J_s$ coincides with $\cO(U(1))$ viewed as a coalgebra,
and that
$\cO(S^2_{q,s})\inc\cO(SU_q(2))$ is a principal
${\cO(SU_q(2))/J_s}$-extension \cite{b-t96,ms99,bm00}.
Next,  the representations \cite{p-p87}
\[
\pi_-(K)e_n=-s^2q^{2n}e_n,\;\;\; \pi_-(L)e_n=\hl_n^-(q,s)e_{n-1},\;\;\;
\hl_n^-(q,s)=s\sqrt{1-(1-s^2)q^{2n}-s^2q^{4n}}
\]
\[
\pi_+(K)e_n=q^{2n}e_n,\;\;\; \pi_+(L)e_n=\hl_n^+(q,s)e_{n-1},\;\;\;
\hl_n^+(q,s)=\sqrt{s^2+(1-s^2)q^{2n}-q^{4n}},
\]
form a needed 1-summable Fredholm module \cite{mnw91},
 and we have our second main
result:
 \bth
For all $\mu\in\Z$, the pairing between the cyclic 0-cocycle
${\rm tr}_\pi$ and the $K_0$-class of $\O(SU_{q}(2))_{\mu,s}$
(Chern number) coincides with
the winding number $\mu$, i.e.,
$\mbox{$
\langle \tr_\pi,[\O(SU_{q}(2))_{\mu,s}]\rangle=\mu.
$}$
\ethe

\section*{Introduction}
\setcounter{equation}{0}

Herein we study two families of noncommutative deformations of the Hopf
fibration $S^3\ra S^2$. The first one  is based on the idea of local
triviality. We can view $S^2$ as the gluing of two discs along their
boundaries, and $S^3$ as the gluing of two solid tori along their boundaries
(a Heegaard splitting of $S^3$).
Then the classical discs can be replaced by Klimek-Lesniewski quantum discs
\cite{kl93}, and subsequently a two-parameter family of noncommutative Hopf
fibrations can be constructed \cite{cm02,hms2}.
The other family originates from the theory of quantum groups. First, one
considers $S^3$ as $SU(2)$ and deforms it into the quantum group $SU_q(2)$
\cite{w-sl87}. Then the classification of $SU_q(2)$-quantum homogeneous spaces
yields a two-parameter family of noncommutative two-spheres \cite{p-p87}. The
latter form the base of the Hopf fibrations of $SU_q(2)$ \cite{bm00}.

Function algebras of total spaces of principal $U(1)$-bundles always decompose
into direct sums of sections of all associated line bundles. The same
phenomenon occurs for both of the aforementioned deformations, i.e., the
coordinate algebras $\cO(S^3_{pq})$ and $\cO(SU_q(2))$ are direct sums of
associated finitely generated projective modules. The aim of this paper is to
prove that these modules are mutually non-isomorphic.

To achieve this, we take an advantage of the K-homology Chern
character \cite{c-a85} and the Chern-Galois character \cite{bh}.
The former produces cyclic cocycles out of bounded
representations of the base algebra, and the latter cyclic cycles out of
finite dimensional corepresentations of the structure coalgebra
that define associated modules. The evaluation of these cocycles
on these cycles gives $K_0$-invariants of the modules. The noncommutative index
formula \cite{c-a85} shows that these invariants are indices of Fredholm
operators, so that they have to be integers. This fact is essentially used
in carrying out the computations.

Throughout this paper we work with unital algebras over a field (complex
numbers in the studied examples) and adopt the standard Hopf-algebraic notation
$m,\hD,\he,S$ for the multiplication, comultiplication, counit and antipode,
respectively. The symbolic notation $\cO(\mbox{quantum space})$ means a
polynomial algebra  defined by generators and relations, and $\Tr$ stands for
the operator trace.

\section{Principal extensions and summable Fredholm modules}
\setcounter{equation}{0}

Let $C$ be a coalgebra and $P$ an algebra and a right $C$-comodule via
$\Delta_R:P\ra P\ot C$. Put
$
B=P^{co C}:=\{b\in P~|~\Delta_R(bp)=b\Delta_R(p),\ \forall p\in P\}.
$
We say that the inclusion $B\inc P$ is a
$C$-extension. A $C$-extension  $B\inc P$ is called {\em principal}
\cite{bh} iff
\begin{enumerate}
\item
the Galois map
$
P\ot_BP\stackrel{can}{\lra} P\ot C,\;p\ot p'\mapsto p\dr(p')
$
is bijective;
\item
there exists a left $B$-linear right $C$-colinear splitting of the
multiplication map $B\ot P\ra P$;
\item
the canonical entwining map $C\ot P\stackrel{\psi}{\ra} P\ot C$, $c\ot p\mapsto
can(can^{-1}(1\ot c)p)$ is bijective;
\item
there is a group-like $e\in C$ such that $\dr(p)=\psi(e\ot p)$, $\forall p\in
P$.
\end{enumerate}
In order to define a strong connection on a principal extension, first we need
to define a bicomodule structure on $P\ot P$.
 The tensor product $P\ot P$ is a
right $C$-comodule via $\hD_R^\otimes:=\id\ot\dr$. Since $\psi$ is bijective,
it is also a left $C$-comodule via $\hD_L^\otimes:=(\psi^{-1}\ci(\id\ot
e))\ot\id$.  The two coactions evidently commute.
Now, let $\pi_B:P\ot P\ra P\ot_B P$ be the canonical surjection.
A linear map $\ell:C\ra P\ot P$ is called a {\em strong connection}
\cite{bh}
 iff it satisfies $can\ci\pi_B\ci\ell=1\ot\id$,
$\hD_R^\otimes\ci\ell=(\ell\ot\id)\ci\hD\;$,
$\hD_L^\otimes\ci\ell=(\id\ot\ell)\ci\hD$,
and
$\ell(e)=1\ot 1$.
If $B\inc P$ is a principal $C$-extension
and $\phi:V_\phi\ra V_\phi\ot C$ is a finite dimensional corepresentation,
then 
the left $B$-module $\hom^C(V_\phi,P)$ of all
colinear maps from  $V_\phi$ to $P$ is finitely generated projective \cite{bh}.
Denote by $\corep_f(C)$ the space of all finite-dimensional corepresentations of
$C$. Then a principal $C$-extension $B\inc P$ yields a map $\phi\mapsto
[\hom^C(V_\phi,P)]$ from $\corep_f(C)$ to $K_0(B)$. The composition of
this map with the Chern character $K_0(B)\ra HC_{even}(B)$ (e.g., see
\cite{l-jl98}) is called the {\em Chern-Galois character} \cite{bh} and is
denoted by $chg$. Explicitly, in degree 0, we have
$chg_0(\phi)=[c_{\phi}^{<2>}c_{\phi}^{<1>}]$. Here $c_{\phi}$ is the character
of $\phi$, i.e., $c_{\phi}=\sum_{i=1}^{\dim V}e_{ii}$,
$\phi(e_j)=\sum_{i=1}^{\dim V}e_{i}\ot e_{ij}$, $\{e_{i}\}$ a basis of
$V_\phi$, and $\ell(c)=:c^{<1>}\ot c^{<2>}$ (summation understood).

Let us now recall the analytic tool of Fredholm modules that are used in the
sequel. A {\em $p$-summable Fredholm module} over a $*$-algebra $B$ can be
viewed as a pair $(\rho_1,\rho_2)$ of bounded $*$-representations of $B$ such
that $
 \Tr|\rho_1(b)-\rho_2(b)|^p<\infty$
for all $b\in B$ (see \cite[p.88]{c-a85} for related details).
The $K$-homology Chern character assigns to finitely summable Fredholm modules
cyclic cocycles \cite{c-a85}. For $p=1$ it takes a particularly simple form,
notably it turns $(\rho_1,\rho_2)$ into a trace (cyclic 0-cocycle) on $B$ via
the formula $\tr_\rho(b)=\Tr(\rho_1(b)-\rho_2(b))$.

\section{A locally trivial quantum Hopf bundle}
\setcounter{equation}{0}
Let us consider the two-parameter family \cite{cm02}  of
$*$-algebras $\O(S^3_{pq})$, $0\leq p,q\leq 1$, generated by $a$ and $b$
satisfying
${a}^*{a} - q{a} {a}^* = 1-q, \;\;\;
{b}^* {b}-p{b} {b}^* = 1-p, \;\;\;
ab=ba,\;\;\; a^*b=ba^*, \;\;\;
(1-{a}{a}^*)(1-{b} {b}^*)= 0.
$
The $*$-subalgebra of $\O(S^3_{pq})$ generated by $ab$ and $bb^*$
can be identified with the $*$-algebra $\O(S^2_{pq})$ generated by $f_0$
and $f_1$ satisfying
\cite{cm00}
$f_0=f_0^*,\;\;\;
f_1^*f_1-qf_1f_1^*=(p-q)f_0+1-p,\;\;\;
f_0f_1-pf_1f_0=(1-p)f_1,\;\;\;
(1-f_0)(f_1f_1^*-f_0)=0.
$
The isomorphism is given by $f_0\mapsto bb^*$ and $f_1\mapsto ab$.
It follows from \cite[Lemma 4.2]{hms2} that $\O(S^2_{pq})\inc \O(S^3_{pq})$
is a principal $\O(U(1))$-extension.  Moreover, we have
\ble[\cite{hms2}]\label{scl1}
Let $u$ be the unitary generator of $\O(U(1))$.
The linear map $\ell:\O(U(1))\ra \O(S^3_{pq})\ot\O(S^3_{pq})$ given
on the basis elements $u^\mu$, $\mu\in\Z$, by the
formulas
\[
\ell(1)=1\ot 1,\;\;\;
\ell(u)=a^*\ot a + qb(1-aa^*)\ot b^*,\;\;\;
\ell(u^*)= b^*\ot b +pa(1-bb^*)\ot a^*,
\]
\[
\ell(u^\mu)=u^{[1]}\ell(u^{\mu-1})u^{[2]},\;\;\;
\ell(u^{*\mu})=u^{*[1]}\ell(u^{*(\mu-1)})u^{*[2]},\;\;\;\mu>0,
\]
defines a strong connection on $\O(S^2_{pq})\inc
\O(S^3_{pq})$.
\ele
The one-dimensional corepresentations of $\O(U(1))$ are labelled by integers.
Explicitly, we have $\phi_\mu(1)=1\ot u^{-\mu}$. Each $\phi_\mu$ yields a
finitely generated projective $\O(S^2_{pq})$-module (quantum line bundle)
$\O(S^3_{pq})_\mu:=\hom^{\O(U(1))}(\C_{\phi_\mu},\O(S^3_{pq}))$.
On the other hand, we have the following irreducible bounded
$*$-representations of $\cO(S^2_{pq})$
\cite[Proposition 19]{cm00}:
\bea
 &&\label{ro1}
\rho_1(f_0)e_k=(1-p^k)e_k,\;\;\;
\rho_1(f_1)e_k=\sqrt{1-p^{k+1}}\:e_{k+1},\;\;\;k\geq 0;
\\ &&\label{ro2}
\rho_2(f_0)e_k=e_k,\;\;\;
\rho_2(f_1)e_k=\sqrt{1-q^{k+1}}\:e_{k+1},\;\;\;k\geq 0.
\eea
Here $\{e_k\}_{k\geq 0}$ is an orthonormal basis of a separable Hilbert space.
Moreover, we have
\ble\label{trlo}
The pair of representations $(\rho_2,\rho_1)$
given by (\ref{ro1})--(\ref{ro2})  yields a 1-summable
Fredholm module over $\O(S^2_{pq})$, so that
$
{\rm tr}_\rho(f):={\rm Tr}(\rho_2(f)-\rho_1(f))
$
defines a trace on $\O(S^2_{pq})$.
\ele

\bth
For all $\mu\in\Z$, the pairing between the cyclic 0-cocycle
${\rm tr}_\rho$ and the $K_0$-class of $\O(S^3_{pq})_\mu$
(Chern number) coincides with the
winding number $\mu$, i.e.,
$
\langle \tr_\rho,[\O(S^3_{pq})_{\mu}]\rangle = \mu.
$
\ethe
{\em Proof outline:}
The pairing of cyclic cohomology and $K$-theory is given by the evaluation
of a cyclic cocycle on the image of the Chern character. In our case (see
Section~1) it gives
$
\langle \tr_\rho,[\O(S^3_{pq})_{\mu}]\rangle =
\tr_\rho(chg_0(\phi_\mu))=\tr_\rho((u^{-\mu})^{<2>}(u^{-\mu})^{<1>}). $
The last expression can be computed explicitly as a function of $p$ or $q$,
depending on whether $\mu$ is positive or negative. Surprisingly, these
functions can be identified with a certain expression appearing in the index
computation carried out in \cite{h-pm00}. Since the latter is proven therein
to be the constant $\mu$, the assertion of the theorem follows.\epf

The trace $\tr_\rho$ computes the Chern numbers of $\O(S^3_{pq})_{\mu}$.
In order to determine the rank of these modules, we can employ any character
of $\O(S^3_{pq})$. Indeed, let $\delta$ be an algebra homomorphism from
$\O(S^3_{pq})$ to $\C$. (See \cite{hms2} for the classification of irreducible
representations of $\O(S^3_{pq})$, including one-dimensional ones.)
Then $\langle
\delta,[\O(S^3_{pq})_\mu]\rangle=\delta((u^{-\mu})^{<2>}(u^{-\mu})^{<1>})=
\delta((u^{-\mu})^{<1>}(u^{-\mu})^{<2>})=\epsilon(u^{-\mu})=1.$
The last two equalities follow from the general properties $m\ci \ell=\epsilon$
(see \cite{bh}) and $\epsilon(\mbox{group-like})=1$, respectively.
The characters always pair integrally with $K_0$ (e.g., see \cite{l-jl98}).
On the other hand, as $\tr_\rho$ comes from a 1-summable Fredholm module, its
pairing with $K_0$ is an index of a Fredholm operator \cite[p.60]{c-a85},
whence also an integer. Therefore, it follows from the linearity of the pairing
that we have a group homomorphism
$
(\delta,\tr_\rho):K_0(\O(S^2_{pq}))\ra\Z\oplus\Z,\;[p]\mapsto
(\langle \delta,[p]\rangle , \langle \tr_\rho,[p]\rangle).
$
The point here is that
for any $\mu\in\Z$ there exists
a rank one projective module with its Chern number equal
to~$\mu$. More formally, we have

\bco\Label{k+}
The image of the positive cone of $K_0(\O(S^2_{pq}))$ under
$(\delta,\tr_\rho):K_0(\O(S^2_{pq}))\ra\Z\times\Z$
contains $\Z_+\times\Z$.
\eco

\section{Hopf fibrations of \boldmath$SU_q(2)$ over
generic Podle\'s spheres} \setcounter{equation}{0}

In this section, we work with the Hopf $*$-algebra $\O(SU_q(2))$,
$q\in ]0,1[$, generated by
elements $\alpha$ and $\gamma$ satisfying \cite{w-sl87}
$\ha\hg=q\hg\ha,\;\;\; \ha\hg^*=q\hg^*\ha,\;\;\; \hg\hg^*=\hg^*\hg,\;\;\;
\ha^*\ha+\hg^*\hg=1,\;\;\; \ha\ha^*+q^2\hg\hg^*=1,
$
and with the Podle\'s spheres \cite{p-p87}. A uniform description of all
Podle\'s spheres  can be obtained by rescaling generators used in \cite{p-p87}.
Then the coordinate $*$-algebras of the quantum spheres $S^2_{q,s}$,
$s\in[0,1]$, are defined by generators $K$ and $L$ satisfying the
relations
$
K=K^*,\;\;\; LK=q^2KL,\;\;\; L^*L+K^2=(1-s^2)K+s^2,\;\;\;
LL^*+q^4K^2=(1-s^2)q^2K+s^2.
$
We can view $\O(S^2_{q,s})$ as a subalgebra of $\O(SU_q(2))$ via the formulas
$
K=s(\hg\ha+\ha^*\hg^*)+(1-s^2)\hg^*\hg,\;\;\;L=s(\ha^2-q{\hg^*}^2)+(1-s^2)\ha\hg^*.
$
Next, let us define the quotient coalgebra $\cO(SU_q(2))/J_s$, where $J_s$ is
the coideal right ideal generated by $K$, $L-s$, $L^*-s$. One can show
that it coincides with $\cO(U(1))$ viewed as a coalgebra, and that
$\cO(S^2_{q,s})=\cO(SU_q(2))^{co\cO(SU_q(2))/J_s}$ \cite{b-t96,ms99}.
Moreover, one can prove that $\cO(S^2_{q,s})\inc\cO(SU_q(2))$
is a principal $\cO(SU_q(2))/J_s$-extension \cite{bm00}.
With the help of  \cite{bh}, the latter follows  from an explicit construction
of a strong connection:
\ble[\cite{bm00}]\label{scl2}
Let $i:\O(U(1))\ra \O(SU_q(2))$ be the linear map defined
on the basis elements $u^\mu$, $\mu\in\Z$, by the
formulas
\[\Label{split}
i(u^{-\mu}):=\left\{\ba{cr}
\prod_{j=0}^{-\mu-1}h_j&\mbox{for } \mu<0\\
1&\mbox{for } \mu=0\\
\prod_{j=0}^{\mu-1}k_j&\mbox{for } \mu>0,
\ea\right.\;\;\;\mbox{(products increase from left to right)}
\]\[\Label{hk}
h_j:=\frac{\ha+q^js(\hg-q\hg^*)+q^{2j}s^2\ha^*}{1+q^{2j}s^2},\;\;\;
k_j:=\frac{\ha^*-q^{-j}s(\hg-q\hg^*)+q^{-2j}s^2\ha}{1+q^{-2j}s^2}.
\]
Then $\ell=(S\ot\id)\ci\hD\ci i$ is a strong connection on
$\O(S^2_{q,s})\inc \O(SU_{q}(2))$.
\ele
Much as before, we associate to every corepresentation $\phi_\mu$
a finitely generated projective left $\O(S^2_{q,s})$-module
$\O(SU_q(2))_{\mu,s}:=\hom^{\cO(SU_q(2))/J_s}(\C_{\phi_\mu},\cO(SU_q(2)))$.
On the other hand, using the representations \cite{p-p87}
\[
\pi_-(K)e_n=-s^2q^{2n}e_n,\;\;\; \pi_-(L)e_n=\hl_n^-(q,s)e_{n-1},\;\;\;
\hl_n^-(q,s)=s\sqrt{1-(1-s^2)q^{2n}-s^2q^{4n}}
\]
\[
\pi_+(K)e_n=q^{2n}e_n,\;\;\; \pi_+(L)e_n=\hl_n^+(q,s)e_{n-1},\;\;\;
\hl_n^+(q,s)=\sqrt{s^2+(1-s^2)q^{2n}-q^{4n}},
\]
one can prove
\ble[\cite{mnw91}]\Label{fred2}
For any $q\in ]0,1[,\; s\in[0,1]$, the pair of representations $(\pi_-,\pi_+)$
yields a 1-summable Fredholm module over $\cO(S^2_{q,s})$,
so that $\tr_\pi:=\Tr\ci(\pi_--\pi_+)$ is a trace on
$\cO(S^2_{q,s})$.
\ele

\bth
For all $\mu\in\Z$, the pairing between the cyclic 0-cocycle
${\rm tr}_\pi$ and the $K_0$-class of $\O(SU_{q}(2))_{\mu,s}$
(Chern number) coincides with
the winding number $\mu$, i.e.,
$
\langle \tr_\pi,[\O(SU_{q}(2))_{\mu,s}]\rangle = \mu.
$
\ethe
{\em Proof outline:}
The proof rests on the following three facts:
$\langle \tr_\pi,[\O(SU_{q}(2))_{\mu,s}]\rangle$ is a rational function of $q$
and $s$,  it is an
integer, and
$\langle \tr_\pi,[\O(SU_{q}(2))_{\mu,s}]\rangle(q,0)=\mu$.
The first claim can be proven with the help of the Chern-Galois
character \cite{bh}, the second follows from
the noncommutative index formula \cite[p.60]{c-a85}, and the third has been
obtained in \cite[Theorem 2.1]{h-pm00}.
Since an integer-valued rational function on a connected set has to be
constant, we can conclude that  $\langle
\tr_\pi,[\O(SU_{q}(2))_{\mu,s}]\rangle(q,s)= \langle
\tr_\pi,[\O(SU_{q}(2))_{\mu,s}]\rangle(q,0)=\mu$.
\epf

\noindent
\bco\Label{k+}
The image of the positive cone of $K_0(\asq)$ under
$(\he,\tr_\pi):K_0(\asq)\ra\Z\times\Z$
contains $\Z_+\times\Z$. (Here $\epsilon$ is the counit of $\O(SU_q(2))$.)
\eco

\footnotesize
{\bf Acknowledgements.}
This work has been partially supported by
a Marie Curie Fellowship
 HPMF-CT-2000-00523 (P.M.H.),
Universit\"at Leipzig (P.M.H., W.S.),
 Deutsche
Forschungsgemeinschaft (R.M.),
Research Grants
Committee of the University of Newcastle and Max-Planck-Institut f\"ur
Mathematik  Leipzig (W.S.). All three authors are
grateful to Mathematisches Forschungsinstitut Oberwolfach for support via its
Research in Pairs programme.

\end{document}